\documentclass[12pt,reqno]{amsart}
\usepackage{amsmath, amsfonts, color, amsthm, graphics, amssymb}
\usepackage[notcite,notref,final]{showkeys}

\numberwithin{equation}{section}
 \theoremstyle{plain}
\newtheorem{theorem}[equation]{Theorem}

\newtheorem{lemma}[equation]{Lemma}
\newtheorem{subclaim}[equation]{Subclaim}
\newtheorem{claim}[equation]{Claim}
\newtheorem{corollary}[equation]{Corollary}
\newtheorem{prop}[equation]{Proposition}

\theoremstyle{definition}
\newtheorem{definition}[equation]{Definition}

\newtheorem{remark}[equation]{Remark}

\newtheorem{example}[equation]{Example}
\newtheorem{examples}[equation]{Examples}

\newcommand{\mf}{\mathfrak}
\newcommand{\Z}{\mathbb{Z}}

\newcommand{\C}{\mathbb{C}}
\newcommand{\N}{\mathbb{N}}

\newcommand{\ten}{\otimes}

\newcommand{\ra}{\rightarrow}
\newcommand{\lra}{\longrightarrow}

\newcommand{\ds}{\displaystyle}
\newcommand{\mc}{\mathcal}

\newcommand{\Aut}{\textrm{Aut}}

\input xy
\xyoption{all}

\begin{document}

\title[Representation theory of $Skly_3$]{Representation theory of three-dimensional Sklyanin algebras}
\author{Chelsea Walton}

\address{Department of Mathematics, University of Washington, Seattle, WA 98195}
\email{notlaw@math.washington.edu}



\bibliographystyle{abbrv}       

\begin{abstract}

We determine the dimensions of the irreducible representations of the Sklyanin algebras with global dimension 3. This contributes to the study of marginal deformations of the N=4 super Yang-Mills theory in four dimensions in supersymmetric string theory. Namely, the classification of such representations is equivalent to determining the vacua of the aforementioned deformed theories. 

We also provide the polynomial identity degree for the Sklyanin algebras that are module finite over their center. The Calabi-Yau geometry of these algebras is also discussed. 

\end{abstract}

\maketitle



\section{Introduction} 

During the last several years, much attention has been paid to supersymmetric gauge theories, and the mathematical results in the work are prompted by these physical principles. 
In particular, we provide representation-theoretic, geometric, and algebraic results about so-called {\it three-dimensional Sklyanin algebras} (Definition \ref{def:Skl3}), all of which are motivated by Leigh-Strassler's seminal work on supersymmetric quantum field theories \cite{LeighStrassler}. Here, it was demonstrated that marginal deformations of the $N=4$ super Yang Mills (SYM) theory in four dimensions yield a class of $N=1$ conformal field theories that are conformally invariant. 
A few years later, Berenstein-Jejjala-Leigh continues this study by considering deformations of the $N=4$ SYM theory in an attempt to classify their moduli space of vacua \cite{BJL}. Although the full classification of vacua was not achieved, the authors introduce a noncommutative geometric framework to study vacua remarkably without the use of $C^{\ast}$-algebras. 

The work in this article begins with the observation that, since an $N=4$ SYM theory can be described in terms of an $N=1$ SYM theory equipped with a superpotential, the problem of classifying vacua reduces to computing matrix-valued solutions to $F$-flatness. As in \cite{BJL}, we aim to understand the holomorphic structure of the moduli space of vacua, so we assume that when given a solution to $F$-flatness that there is also a solution to the $D$-term constraints. Hence, determining vacua is algebraically equivalent to determining finite-dimensional matrix-valued solutions to the system of equations $\{\partial_{\text cyc}\Phi = 0\}$, where $\partial_{cyc}\Phi$ is the set of cyclic derivatives of the superpotential $\Phi$ derived from the $N=1$ SYM theory mentioned above. 

To see how three-dimensional Sklyanin algebras arise in this study, we consider marginal deformations of the $N=4$ SYM theory studied in \cite{BJL}. According to \cite{BJL}, these deformations can be presented as an $N=1$ SYM theory with homogeneous superpotential:
$$\Phi_{marg} = axyz + byxz + \frac{c}{3}(x^3 + y^3 + z^3),$$
where $a$, $b$, and $c$ are scalars. Since finite-dimensional matrix solutions to $\{\partial_{\text cyc}\Phi_{marg} = 0\}$ are desired, we are interested in computing finite-dimensional irreducible representations (irreps) of an algebra whose ideal of relations is $(\partial_{\text cyc}\Phi_{marg})$. More precisely, we have that the classification of the moduli space of vacua of marginal deformations of the $N=4$ SYM theory in four dimensions boils down to classifying irreps of the following {\it superpotential algebras} (Definition \ref{def:suplalg}).

\begin{definition} \label{def:Skl3} Let $k$ be an algebraically closed field of characteristic not equal to 2 or 3. The {\it three-dimensional Sklyanin algebras}, denoted by $S(a,b,c)$ or $Skly_3$, are generated by three noncommuting variables $x$, $y$, $z$, subject to three relations:
\begin{equation} \label{eq:S(abc)}
\begin{array}{c}ayz+ bzy+ cx^2 ~=~ azx+ bxz+cy^2 ~=~ axy+ byx+ cz^2=0. \end{array}
\end{equation}
Here, we require that:\\ 
\indent (1) $[a:b:c] \in \mathbb{P}^2_k \setminus \mf{D}$ where 
$$\mf{D} = \{[0:0:1], [0:1:0], [1:0:0]\} ~\cup~ \{[a:b:c] ~|~ a^3=b^3=c^3=1\};$$

(2) $abc \neq 0$ with $(3abc)^3 \neq (a^3+b^3+c^3)^3$.

\end{definition}

Results on the dimensions of irreps of $Skly_3 = S(a,b,c)$ are summarized in Theorem \ref{thm:Chap4sum} below, yet first we  must recall for the reader the noncommutative geometric data associated to these structures. In the sense of noncommutative projective algebraic geometry, the algebras $Skly_3$ come equipped with geometric data ($E$, $\sigma$).  Here,  $E = E_{abc}$ is an elliptic curve in projective space $\mathbb{P}^2$, and $\sigma= \sigma_{abc}$ is an automorphism on $E$. The algebras $S=Skly_3$ also contain a central element $g$ so that $S/Sg \cong B$, where $B$ is a {\it twisted homogeneous coordinate ring} (Definition \ref{def:thcr}, Theorem \ref{thm:ATVSkly3}). Details about the noncommutative geometry of the rings $S$ and $B$ are provided in $\S$2.3.
Now, we state the main result on the dimensions of irreducible representations of $Skly_3$.

\begin{theorem} \label{thm:Chap4sum} (Lemma \ref{lem:inf.order}, Theorem \ref{thm:fin.ord.g-tfree}, Theorem \ref{thm:fin.ord.g-tors})\\
(i) When $|\sigma| = \infty$, the only finite-dimensional irrep of $Skly_3$ is the trivial representation. 

\noindent (ii) Assume that $|\sigma|<\infty$. Let $\psi$ be a nontrivial irrep of $Skly_3$. 
\begin{enumerate}
\item[(a)] If $\psi$ is $g$-torsionfree as defined in $\S$2.1, then $\dim_k \psi = |\sigma|$ when  $(3,|\sigma|) =1$; and $|\sigma|/3 \leq \dim_k \psi \leq |\sigma|$ when $(3,|\sigma|)\neq 1$. 
\item[(b)] If $\psi$ is $g$-torsion, then $\psi$ is a nontrivial irrep of the twisted homogeneous coordinate ring $B$, and $\dim_k \psi = |\sigma|$.
\end{enumerate}
\end{theorem}
 
Here, the automorphism $\sigma$ is generically of infinite order; thus, the representation theory of $Skly_3$ is generically trivial. On the other hand, higher-dimensional representations of $Skly_3$, or higher-dimensional matrix-valued solutions to $F$-flatness, only arise when $|\sigma| < \infty$. The dimension of an irrep is generically equal to $|\sigma|$ in this case, and surprisingly there exist irreps of $Skly_3$ of intermediate dimension.

\begin{prop} (Proposition \ref{prop:Simp2})
Given parameters $(a,b,c) = (1,-1,-1)$, the automorphism $\sigma_{1,-1,-1}$ has order 6 and there exist 2-dimensional irreps of the three-dimensional Sklyanin algebra $S(1,-1,-1)$. 
\end{prop}

\noindent The occurrence of the automorphism $\sigma$ with finite order is also key to understanding the Calabi-Yau (CY) geometry of deformed $N=4$ SYM theories as we will now see.

The geometric results in this article are prompted by Maldacena's AdS/CFT correspondence \cite{Maldacena}. Algebras that are module-finite over their center play a crucial role in this result as algebras with large centers yield $AdS_5 \times X_5$ without flux. (Here, $X_5$ is a five-dimensional Sasaki-Einstein manifold.) Furthermore in this case, the centers of such algebras are typically CY cones over $X_5$: a cone due to the grading of the algebra, and a CY variety due to conformal symmetry.

In light of the discussion above, naturally we ask which of the three-dimensional Sklyanin algebras are module-finite over their center, or which of these algebras potentially have Calabi-Yau geometric structure. To answer this question, we note that $S(a,b,c)$ is module-finite over its center if and only if $\sigma = \sigma_{abc}$ has finite order \cite[Theorem 7.1]{ATV1}. Thus, we aim to classify parameters $(a,b,c)$ for which $|\sigma_{abc}|=n <\infty$. This has been achieved for small values of $n$.

\begin{prop} (Proposition \ref{prop:sigmafinite}) The parameters $(a,b,c)$ for which $|\sigma_{abc}|=n$ has been determined for $n=1,\dots,6$.
\end{prop}

\noindent For some $n$, we have a 1-parameter family of triples $(a,b,c)$ so that $|\sigma_{abc}|=n$. In $\S$5, we illustrate how some degeneration limits of these 1-parameter families arise in the analysis of orbifolds with discrete torsion \cite[$\S$4.1]{BJL}. The centers $Z$ of $S(a,b,c)$ at some of these degeneration limits are also presented, and for such $Z$, we have that Spec($Z$) is the affine toric CY three-fold $\mathbb{C}^3/(\Z_n \times \Z_n)$  (Proposition \ref{prop:centerdeg}, Corollary \ref{cor:CY}).

Finally, the algebraic results about $S=Skly_3$ also pertain to the case when $S$ is module-finite over its center, or equivalently when $S$ satisfies a {\it polynomial identity} (is {\it PI})   (Definition~\ref{def:PI}). Note that we have that the twisted homogeneous coordinate ring $B$ of Theorem~\ref{thm:Chap4sum} is PI if and only if $S$ is PI. In the result below, we determine the {\it PI degree} (Definition \ref{def:PIdeg}) of both $S$ and $B$, which is simply a measure of the rings' noncommutativity.

\begin{prop} \label{prop:PIdegintro} (Corollary \ref{cor:PIdegB}) When the automorphism $\sigma$ associated to $S(a,b,c)$ satisfies $|\sigma| < \infty$, then both the PI degree of $S(a,b,c)$ and  of the corresponding twisted homogeneous coordinate ring $B$ are equal to $|\sigma|$.
\end{prop}

\noindent This result also verifies a conjecture of Artin and Schelter (Corollary \ref{cor:conjAS}). Information about PI rings and their associated geometry can be found in $\S$2.2.

Further results of this work are also discussed in $\S$5. For instance, we can consider the {\it relevant} deformations of the SYM theory in the sense of \cite{BJL}. This requires the classification of irreducible representations of {\it deformed Sklyanin algebras}; see $\S$5.2. Moreover in $\S$5.3, we analyze the representation theory of twelve algebras that were omitted from the family of algebras $S(a,b,c)$ in Definition~\ref{def:Skl3}, so-called {\it degenerate Sklyanin algebras}. 

\medskip

\noindent {\bf Acknowledgments.} Most of the results in this article were part of my Ph.D. thesis at the University of Michigan, and I thank my adviser Toby Stafford for his guidance on this project. Moreover, I thank  David Berenstein, Kenneth Chan, David Morrison, Michael Wemyss, and Leo Pando Zayas for supplying many insightful comments.  I am especially grateful to the referee for their extensive feedback.  Much of this work was expanded and edited during the Noncommutative Geometry and $D$-branes workshop at the Simons Center for Geometry and Physics in  December 2011; I thank its participants for many stimulating talks and discussions. I have been supported by the National Science Foundation: grants DMS-0555750, 0502170, 1102548.


\section{Background material}
\label{sec:background}

Here, we provide the background material of this article, which pertains to: (1) the interaction between physics and the representation theory of superpotential algebras; (2)  polynomial identity (PI) rings; and (3) the field of noncommutative projective algebraic geometry. The reader is referred to \cite{McConnellRobson} and \cite{StaffordVandenbergh} for a thorough introduction to the latter two areas, respectively.

\subsection{From physics to representation theory}

As mentioned in the introduction, the study of SYM theories equipped with a superpotential $\Phi$ prompts the analysis of the representation theory of a superpotential algebra. Namely, to satisfy the F-term constraints of the SYM theory, one must determine irreducible matrix solutions to cyclic derivatives of $\Phi$. This task is equivalent to classifying irreducible representations of the corresponding superpotential algebra. We define these terms as follows.

Given an algebra $A$, an {\it $n$-dimensional representation} of $A$ is a ring homomorphism from $A$ to the ring of $n \times n$ matrices, $\psi: A \ra Mat_n(k)$. 
We say that $\psi$ is {\it irreducible}, or is an {\it irrep}, if there does not exist a proper subrepresentation of $\psi$. Moreover, the set of irreps of $A$ is in bijective correspondence to simple left $A$-modules. We call $\psi$ {\it g-torsion} (or {\it g-torsionfree}, respectively) if the corresponding module $M$ satisfies $gM=0$ (or $gm\neq0$ for all $0 \neq m \in M$, respectively).

Now the algebras under consideration in this work, namely superpotential algebras, are defined as follows. 

\begin{definition}  \label{def:suplalg} Let $V$ be a $k$-vector space with basis $x_1, \dots, x_n$ and let $F = T(V) = k\{x_1,\dots,x_n\}$ be the corresponding free algebra. The commutator quotient space $F_{cyc} = F/[F,F]$ is a $k$-vector space with the natural basis formed by cyclic words in the alphabet $x_1, \dots, x_n$. 
Elements of $F_{cyc}$ are referred to as {\it superpotentials} (or as {\it potentials} in some articles). 
Let $$\Phi = \ds \sum_{\{i_1,i_2,\dots,i_r\} \subset I} x_{i_1} x_{i_2} \cdots x_{i_r} \in F_{cyc}$$ for some indexing set $I$. For each $j = 1, \dots,n$, one defines $\partial_j \Phi \in F$ as the corresponding partial derivative of $\Phi$ given by:
$$\partial_j \Phi = \sum_{\{s | i_s = j\}} x_{i_s + 1} x_{i_s +2} \cdots x_{i_r} x_{i_1} x_{i_2} \cdots x_{i_s -1} ~~ \in F.$$
The elements $\partial_j \Phi$ are called {\it cyclic derivatives} of $\Phi$, and
the algebra $F / \left( \partial_j \Phi\right)_{j=1,\dots,n}$ is called a {\it superpotential algebra}.
\end{definition}

\begin{example} Let $F = k\{x,y,z\}$ with superpotential $\Phi = xyz-xzy$. We have that
$\partial_x \Phi = yz -zy$,
$\partial_y \Phi = zx -xz$, and 
$\partial_z \Phi = xy -yx$.
Hence, the polynomial algebra $k[x,y,z]$ arises as the superpotential algebra $F/(\partial_{cyc} \Phi)$.
\end{example}

\begin{example} \cite{BJL}
Consider $F = k\{x,y,z\}$ with the superpotential $$\Phi_{marg} = axyz + byxz + \frac{c}{3}(x^3+y^3+z^3).$$ We have that $\partial_x \Phi_{marg}$, $\partial_y \Phi_{marg}$, and $\partial_z \Phi_{marg}$ are precisely the relations of the three-dimensional Sklyanin algebra (Definition \ref{def:Skl3}). Hence, for generic parameters $(a,b,c)$, $Skly_3$ arises as the superpotential algebra $F/(\partial_{cyc} \Phi_{marg})$.
\end{example}

Many other algebras arise as superpotential algebras, including quantum planes \cite[Chapter 1]{GoodearlWarfield}, the conifold algebra \cite[Example 1]{BLS}, and graded Calabi-Yau algebras of dimension 3 \cite{Bocklandt:CY3}.

\subsection{Polynomial identity rings}

We now discuss a class of noncommutative algebras which have a rich representation theory, algebras that are module-finite over their centers. In fact, we will see that the geometry of the centers of these algebras controls the behavior of their irreducible representations. We apply this discussion specifically to a certain subclass of Sklyanin algebras, those with automorphism $\sigma$ of finite order (see $\S\S$4,5). Rings that are module-finite over their center  are known to satisfy a {\it polynomial identity}, so we state the following results in the language of PI theory. We review the definition, properties, and examples of PI rings as presented in \cite[Chapter 13]{McConnellRobson} as follows.

\begin{definition} \label{def:PI} \cite[13.1.2]{McConnellRobson} A {\it polynomial identity (PI) ring} is a ring $A$ for which there exists a monic multilinear polynomial $f \in \Z\{x_1, \dots, x_n\}$ so that $f(a_1,\dots, a_n)  =0$ for all $a_i \in A$. The minimal degree of such a polynomial is referred as the {\it minimal degree} of $A$. 
\end{definition}

Any ring that is module-finite over its center is PI, and any subring or homomorphic image of a PI ring is PI \cite[13.0]{McConnellRobson}.

\begin{examples}\label{ex:PI1} (1) Commutative rings $R$ are PI as the elements satisfy the polynomial identity $f(x,y) = xy-yx$ for all $x, y \in R$.\\
(2) The ring of $n \times n$ matrices over a commutative ring $R$, Mat$_n(R)$, is PI.\\
(3) The ring $A = k\{x,y\}/(xy+yx)$ is PI as $A$ is a finite module over its center $Z(A) = k[x^2,y^2]$.\\
(4)  If $A$ is PI, then the matrix ring Mat$_n(A)$ is also PI for all $n \geq 1$.
\end{examples}

The notion of the {\it PI degree} of a noncommutative ring can be defined for the class of {\it prime} (\cite[$\S$3.1]{GoodearlWarfield}) PI rings via Posner's Theorem below. Roughly speaking, this degree measures the noncommutativity of a ring: it is equal to 1 for commutative rings and is greater than 1 for noncommutative rings. We refer to \cite[Chapter 6]{GoodearlWarfield} for the definition of a {\it quotient field} and of a {\it Goldie quotient ring}.

\begin{theorem} (Posner) \label{thm:Posner} \cite[13.6.5]{McConnellRobson} Let $A$ be a prime PI ring with center $Z$ with minimal degree $2p$. For $\mc{S} = Z \setminus \{0\}$, let $Q = A\mc{S}^{-1}$ denote the Goldie quotient ring of $A$.  Moreover, let $F = Z\mc{S}^{-1}$ denote the quotient field of $Z$. Then, $\dim_F Q = p^2$. 
\end{theorem}

\begin{definition} \label{def:PIdeg}
Consider the notation of Theorem \ref{thm:Posner}. We say that the {\it PI degree of a prime PI ring} $A$ is equal to $(\dim_F Q)^{1/2}$.
\end{definition}

\begin{examples} \label{ex:PI2} (1) Let $R$ be a commutative ring. Then, $R$ has minimal degree 2,  has PI degree 1, and $R$ is rank 1 over its center.\\
(2) Let $R$ be a commutative ring. Then, Mat$_n(R)$ has minimal degree $2n$ and PI degree $n$, and Mat$_n(R)$  is rank $n^2$ over its center.\\
(3) The ring $A = k\{x,y\}/(xy+yx)$ has minimal degree $4$ and PI degree $2$, and $A$ is rank $4$ over its center.
\end{examples}

Motivating the geometry of PI rings, first we note that irreps of finitely generated PI $k$-algebras $A$ are  finite-dimensional \cite[13.10.3]{McConnellRobson}. Secondly, if $\psi$ is an irrep of $A$, then by Schur's lemma \cite[0.1.9]{McConnellRobson} there exists a unique maximal ideal $\mf{m}$ in the maximal spectrum {\tt maxSpec}($Z(A)$) corresponding to $\psi$. Conversely, the induced map from the set of isomorphism classes of $A$, denoted {\tt Irrep}($A$), to {\tt maxSpec}($Z(A)$) is surjective with finite fibers. Thus, the affine variety {\tt maxSpec}($Z(A)$) plays a crucial role in understanding representation theory of $A$. Refer to Example \ref{ex:GeomPI} below for an illustration of this discussion.

To understand the dimensions of irreps of prime PI algebras geometrically, we employ the proposition below. We say that a prime PI algebra $A$ is {\it Azumaya} if PIdeg($A/\mf{p}$) = PIdeg($A$) for all prime ideals $\mf{p}$ of $A$; see \cite[$\S$13.7]{McConnellRobson} for more details.

\begin{prop} \label{prop:dimnsimple} \cite[Proposition 3.1]{BrownGoodearl} \cite[$\S$III.1]{book:BrownGoodearl}
Let $A$ be a prime noetherian finitely generated $k$-algebra that is module finite over its center $Z$. The following statements hold.

\noindent (a) The maximum $k$-vector space dimension of an irrep of $A$ is PIdeg($A$).

\noindent (b) Let $\psi$ be an irrep of $A$, let $P$ denote ker($\psi$), and let $\mf{m} = P \cap Z$. Then $\dim_k \psi =$ PIdeg($A$) if and only if $A_{\mf{m}}$ is an Azumaya over $Z_{\mf{m}}$. Here, $A_{\mf m}$ is the localization $A \ten_Z Z_{\mf m}$.

\noindent (c) The ideal $\mf{m}$ in (b) is referred to as an {\it Azumaya point}, and the set of Azumaya points is an open, Zariski dense subset of {\tt maxSpec}($Z$).
\end{prop}

In short, the generic irreps of a PI algebra $A$ have maximal dimension, a quantity which is equal to the PI degree of $A$. 

\begin{example} \label{ex:GeomPI} Consider the algebra $A = k\{x,y\}/(xy+yx)$ from Examples~\ref{ex:PI1}(3)~and~\ref{ex:PI2}(3). We have that {\tt maxSpec}$(Z(A))$ is the affine space $\mathbb{A}^2$. By Proposition \ref{prop:dimnsimple}, the maximum dimension of the irreps of $A$ equals the PI degree of $A$, which is equal to 2. 

We can list the nontrivial the 1- and 2-dimensional irreps of $A$, $\{\psi_1\}$ and $\{\psi_2\}$ respectively, as follows. The irreps
$\psi_1: A \ra \text{Mat}_1(k)$ are given by $\psi_1(x) = \alpha$ and $\psi_1(y) = \beta$ where $(\alpha, \beta) = (\alpha,0)$ or $(0,\beta)$. Moreover, the irreps
$\psi_2: A \ra \text{Mat}_2(k)$ are  given by 
\begin{center}
$\psi_2(x)$ = {\footnotesize$\left(\begin{array}{cc}\alpha & 0\\0& -\alpha \end{array}\right)$} ~and~ $\psi_2(y)$ = {\footnotesize$\left(\begin{array}{cc}0 & \beta\\ -\beta& 0 \end{array}\right)$}
\end{center}
where $(\alpha, \beta) \neq (0,0)$.

Thus, the origin of $\mathbb{A}^2$={\tt maxSpec}$(Z(A))$ corresponds to the trivial irrep of $A$. The axes of $\mathbb{A}^2$, save the origin, correspond to the 1-dimensional irreps of $A$. Finally, the bulk of the irreps of $A$, the ones of dimension 2, correspond to the points $(\alpha, \beta) \in \mathbb{A}^2$ where $\alpha,\beta \neq 0$. 
\end{example}

\subsection{Noncommutative projective algebraic geometry}

Here, we provide a brief description of Artin-Tate-van den Bergh's projective geometric approach to understanding noncommutative graded rings; see \cite{ATV1} for details.
We consider algebras $A$ that are {\it connected graded}, meaning that $A=\bigoplus_{i \in \N} A_i$ is $\N$-graded with $A_0=k$. Furthermore, we require that $\dim_k A_i<\infty$ for all $i \geq 0$; in other words $A$ is {\it locally finite}. We also consider the following module categories of $A$:

\smallskip

$A${\it-gr}, the category of $\Z$-graded left noetherian $A$-modules $M = \bigoplus_i M_i$ with degree preserving homomorphisms;

$A${\it-qgr}, the quotient category $A$-gr/tors($A$),  where tors($A$) is the category of finite dimensional graded $A$-modules. 
The objects of $A$-qgr are those of $A$-gr; here $M \cong N$ if there exists an integer $n$ for which $\bigoplus_{i \geq n} M_i \cong \bigoplus_{i \geq n} N_i$.

\smallskip

 Given a connected graded, locally finite algebra $A$, we build a projective geometric object $X$ corresponding to $A$, which will be analogous to the projective scheme Proj($A$) in classical projective algebraic geometry. The closed points of $X$ for $A$ are interpreted as $A$-modules, as two-sided ideals do not generalize properly in the noncommutative setting. More explicitly, a noncommutative point is associated to an {\it $A$-point module}, a cyclic graded left $A$-module $M= \bigoplus_{i \geq 0} M_i$ where $\dim_k{M_i}=1$ for all $i \geq 0$. This module behaves like a homogeneous coordinate ring of a closed point in the commutative setting.
The parameterization of isomorphism classes of $A$-point modules (if it exists as a projective scheme) is called the {\it point scheme} $X$ of $A$, and this is the geometric object that we associate to the noncommutative graded ring $A$. In fact, if $R$ is a commutative graded ring generated in degree 1, then the point scheme of $R$ is Proj($R$). Here are some examples of the discussion above.

\begin{example} \label{ex:ptscheme}
(1) The point scheme of the commutative polynomial ring $R=k[x,y]$ is the projective line $\mathbb{P}^1$ = Proj($R$). Closed points of $\mathbb{P}^1$ are denoted by $[\alpha:\beta]$ for $(\alpha,\beta) \in k^2 \setminus (0,0)$, and are in correspondence with the set of $R$-modules: $\{R/(\alpha x - \beta y) ~|~ [\alpha: \beta] \in \mathbb{P}^1\}.$

\smallskip

\noindent (2) Consider a noncommutative analogue of $k[x,y]$, the  ring $R_q = k\{x,y\}/(yx-qxy)$ for $q \neq 0$, not a root of unity. The set of point modules of $R_q$ is 
$\left\{ R_q / R_q(\alpha x - \beta y) ~|~ [\alpha:\beta] \in\mathbb{P}^1 \right\}.$
Thus, the point scheme of $R_q$ is again $\mathbb{P}^1$.
\end{example}

Returning to algebra, we build a (noncommutative) graded ring $B$ corresponding to the point scheme $X$ of $A$, which often we use to study the ring-theoretic behavior of $A$.

\begin{definition} \label{def:thcr} \cite[$\S$6]{ATV1} Given a projective scheme $X$, let $\mc{L}$ be an invertible sheaf on $X$, and let $\sigma \in \Aut X$. The {\it twisted homogeneous coordinate ring} $B=B(X,\mc{L},\sigma)$ of $X$ with respect to  $\mc{L}$ and $\sigma$ is an $\N$-graded ring: 
$B=\bigoplus_{d \in \N}  H^0(X, \mc{L}_d)$,
 where $\mc{L}_0 = \mc{O}_X$, $\mc{L}_1 = \mc{L}$, and $\mc{L}_d = \mc{L} \ten_{\mc{O}_X} \mc{L}^{\sigma} \ten_{\mc{O}_X} \dots \ten_{\mc{O}_X} \mc{L}^{\sigma^{d-1}}$ for $d \geq 2$. Multiplication is defined by taking global sections of the isomorphism $\mc{L}_d \otimes_{\mc{O}_X} \mc{L}_e^{\sigma^d} \cong \mc{L}_{d+e}$.
\end{definition}

When $\sigma = \text{id}_X$, then we get the commutative {\it section ring} $B(X,\mc{L})$ of $X$ with respect to $\mc{L}$. Therefore, point schemes and twisted homogeneous coordinate rings are respectively genuine noncommutative analogues of projective schemes, Proj, and section rings, $B(X,\mc{L})$, in classical algebraic geometry.

\begin{example} \label{ex:thcr}
Continuing with Example \ref{ex:ptscheme}, let $X = \mathbb{P}^1$ and $\mc{L}=\mc{O}_{\mathbb{P}^1}(1)$. Then for $\sigma$=id$|_{\mathbb{P}^1}$, we get that  $B(X,\mc{L},\text{id}|_{\mathbb{P}^1}) = k[x,y]$. Whereas for $\sigma_q$([$u$:$v$]) = [$qu$:$v$], we have that $B(X,\mc{L},\sigma_q) = R_q$.
\end{example}

The significance of  twisted homogeneous coordinate rings is displayed in the next result.

\begin{theorem} \label{thm:ATVSkly3} \cite{ATV1} \cite[(10.17)]{ArtinSchelter} (i) The point scheme of the three-dimensional Sklyanin algebra $S = S(a,b,c)$ is isomorphic to:
\begin{equation} \label{eq:E(abc)}
E=E_{abc}\index{E_{abc}}: \mathbb{V}\left((a^3+b^3+c^3)xyz - (abc)(x^3+y^3+z^3)\right) ~\overset{i}{\subset}~ \mathbb{P}^2.
\end{equation}
Here,  $E$ is a smooth elliptic curve as $abc \neq 0$ and $(3abc)^3 \neq (a^3+b^3+c^3)^3$; see Definition~\ref{def:Skl3}. The automorphism $\sigma = \sigma_{abc}$ of $E$ is induced by the shift functor on point modules. Moreover, there is a ring surjection from $S(a,b,c)$ to the twisted homogeneous coordinate ring $B(E,i^{\ast}\mc{O}_{\mathbb{P}^2}(1),\sigma)$. 
Even in the case when $abc = 0$, the kernel of the surjection $S \twoheadrightarrow B$ generated by a degree 3 central, regular element given by:
\begin{equation} \label{eq:g}
g = c(a^3-c^3)x^3 + a(b^3-c^3)xyz + b(c^3-a^3)yxz  +  c(c^3-b^3)y^3.
\end{equation}
\end{theorem}

One consequence of Theorem \ref{thm:ATVSkly3} is that three-dimensional Sklyanin algebras are noetherian domains with the Hilbert series of that of a polynomial ring in three variables. Such results were unestablished before the development of noncommutative projective algebraic geometry; see \cite{ATV1} for details.


\section{Irreducible finite dimensional representations of $Skly_3$}

The aim of this section is to establish Theorem \ref{thm:Chap4sum} and Proposition \ref{prop:PIdegintro}. In other words, we will determine the dimensions of irreducible finite-dimensional representations of the three-dimensional Sklyanin algebras $Skly_3$ and of corresponding twisted homogeneous coordinate rings $B$. Moreover, we also compute the PI degrees (Definition \ref{def:PIdeg}) of both $Skly_3$ and $B$. We sometimes denote $Skly_3$ by $S$ or $S(a,b,c)$, and we keep the notation from Theorem \ref{thm:ATVSkly3}.

\medskip

\noindent {\it Notation.} Given an algebra $A$, let {\tt Irrep}$_{< \infty} (A)$ be the set of isomorphism classes of finite-dimensional irreps of $A$, see $\S$2.1. Let {\tt Irrep}$_m (A)$ be the ones of dimension $m$. 

\medskip

The results of this section are given as follows. First, we consider the order of the automorphism $\sigma$, and investigate {\tt Irrep}$_{< \infty}S$ when $|\sigma| = \infty$ (Lemma \ref{lem:inf.order}). For the $|\sigma| < \infty$ case, we also take into consideration the subcases where finite-dimensional irreps of $S$ are either $g$-torsionfree or $g$-torsion. These results are reported in Theorem~\ref{thm:fin.ord.g-tfree} and Theorem~\ref{thm:fin.ord.g-tors}, respectively.

We describe {\tt Irrep}$_{< \infty} S$ in the case that $|\sigma| = \infty$. First we require the following lemma. We use the notion of {\it Gelfand-Kirillov dimension}, a growth measure of algebras and modules in terms of a generating set; refer to \cite{KrauseLenagan} for further details.

\begin{lemma} \label{lem:GKdimA/P} Let $A$ be a finitely generated, locally finite, connected graded $k$-algebra. Take $\psi \in$ {\tt Irrep}$_{< \infty} A$ and let $P$ be the largest graded ideal contained in ker$(\psi)$. Then, the Gelfand-Kirillov (GK) dimension of $A/P$ is equal to 0 or 1.
Moreover, GKdim($A/P$)=0 if and only if $\psi$ is the trivial representation. 
\end{lemma}

\begin{proof}
Let $\overline{A} := A/P$. Let $M$ be the simple left $A$-module corresponding to $\psi \in$ {\tt Irrep}$_{< \infty}(\overline{A})$. Note that ann$_{\overline{A}}(M)$ contains no homogeneous elements. Now suppose by way of contradiction that GKdim($\overline{A}$)$>1$. Then $\lim_{i \ra \infty} \dim_k \overline{A}_i = \infty$. Since $\overline{A}/\text{ann}_{\overline{A}}(M)$ is finite-dimensional by the Density Theorem \cite[Theorem 0.3.6]{McConnellRobson}, we get that ann$_{\overline{A}}M \cap \overline{A}_i \neq 0$ for all $i >> 0$. Hence there exists a homogeneous element in ann$_{\overline{A}}(M)$, which contradicts the maximality of $P$. Now GKdim($\overline{A}$) = 0 or 1.

If $M = A/A_+$, then $P = A_+$, and GKdim($A/P$)=0 as $A/P$ is finite-dimensional. Conversely, GKdim($A/P$) = 0 implies that $\overline{A} = A/P$ is finite-dimensional and connected graded. Since $\bigcap_{n \in \N} \left(\overline{A}_+\right)^n = 0$, we get that $\left(\overline{A}_+\right)^n = 0$ for some $n \in N$. Hence $\overline{A}_+ = 0$ by primality, and so $P = A_+$ with $M$ being the trivial module. 
\end{proof}

\begin{lemma} \label{lem:inf.order}
When $|\sigma| = \infty$, the set $\text{{\tt Irrep}}_{< \infty}(S)$ solely consists of the trivial representation.
\end{lemma}

\begin{proof} 
According to \cite[Lemma 4.1]{SmithStaniszkis}, finite dimensional irreps of $S$ are quotients of irreducible objects in the category $S$-qgr, see $\S$2.3. For $|\sigma| = \infty$, the set of nontrivial, $g$-torsionfree quotients of $S$-qgr is empty due to \cite[Propositions~7.5,~7.9]{ATV2}. 

On the other hand, the set of $g$-torsion irreducible objects of $S$-qgr equals the set of irreducible objects of $B$-qgr, where $B$ is the twisted homogeneous coordinate ring  $B(E,\mc{L},\sigma)$ (Theorem \ref{thm:ATVSkly3}). 
Take $\psi \in$ {\tt Irrep}$_{< \infty}B$ and let $P$ be the largest graded ideal contained in ker($\psi$). By Lemma \ref{lem:GKdimA/P}, we have that GKdim($B/P$)$\leq 1$. If GKdim($B/P$) equals 1, then the Krull dimension (\cite[Chapter 6]{McConnellRobson}) Kdim($B/P$) also equals 1. This is a contradiction as $|\sigma| = \infty$, and $B$ is projectively simple in this case \cite{RRZ}. Thus, GKdim($B/P$)=0, and again by Lemma \ref{lem:GKdimA/P} we know that $\psi$ is the trivial representation.
\end{proof}

We now determine the dimensions of finite-dimensional irreps of $S$ in the case where $|\sigma|<\infty$. We require the following preliminary result, which is mentioned in \cite[$\S$3]{LeBruyn:centsing} with details omitted. The full proof is provided in \cite[Lemma IV.14]{Walton:thesis}. 

\medskip

\noindent {\it Notation.} 
Put $\Lambda:= S[g^{-1}]$, where $g$ is defined in Theorem \ref{thm:ATVSkly3}, and $\Lambda_0$ its degree 0 component. 
By \cite[Theorem 7.1]{ATV2}, $S$ is module finite over its center precisely when $|\sigma|=n < \infty$. Thus in this case, $S$ is PI and we denote its PI degree by $p$ (Definitions \ref{def:PI} and \ref{def:PIdeg}). 

\smallskip

\begin{prop} \label{prop:pideglam0} \cite{LeBruyn:centsing} When $|\sigma|<\infty$,  PIdeg$(\Lambda_0) = p/$gcd$(3,|\sigma|)$. \qed
\end{prop}

\begin{corollary} \label{cor:PIdegS}
When $\sigma = \sigma_{abc}$ has finite order,  the PI degree of $Skly_3=S(a,b,c)$ is equal to $|\sigma_{abc}|$.
\end{corollary}

\begin{proof} 
Let $s$ denote the smallest integer such that $\sigma^s$ fixes $[\mc{L}]$ in Pic$E$.  We know by \cite[Theorem 7.3]{ATV2} that the ring $\Lambda_0$ is Azumaya \cite[Definition 13.7.6]{McConnellRobson} and $s$ is its PI degree. 
Now suppose that $(3,|\sigma|)=1$. Then $s$ = $|\sigma|$ by \cite[$\S$5]{Artin}. Therefore with Proposition \ref{prop:pideglam0}, we have that
$\text{PIdeg }S = \text{PIdeg }\Lambda_0 = s = |\sigma|.$

On the other hand, suppose that $|\sigma|$ is divisible by 3. Now $s$ is also the order of the automorphism $\eta$ introduced in \cite[$\S$5]{ATV2}, due to \cite[Theorem 7.3]{ATV2} (or more explicitly by \cite[Lemma 5.5.5(i)]{ArtinDeJong}). Since $\mc{L}$ has degree 3, we know that $\eta$ is $\sigma^3$ \cite[Lemma 5.3.6]{ArtinDeJong}. Therefore,
$$\text{PIdeg }S ~=~ 3\cdot \text{PIdeg }\Lambda_0
 ~=~ 3|\eta| ~=~ 3|\sigma^3| ~=~ 3 \cdot \frac{|\sigma|}{(3,|\sigma|)} ~=~  |\sigma|.$$
 
\vspace{-.2in} 
\end{proof}

Now we consider two subcases of the classification of dimensions of irreps of $S$ for $|\sigma|<\infty$: first consisting of $g$-torsionfree irreps in Theorem \ref{thm:fin.ord.g-tfree}, and secondly consisting of $g$-torsion irreps in Theorem \ref{thm:fin.ord.g-tors}.

\begin{theorem} \label{thm:fin.ord.g-tfree} Let $|\sigma|<\infty$, and let $\psi$ be a $g$-torsionfree finite-dimensional irrep of $S$. Then, $\dim_{k} \psi=|\sigma|$ for $(3,|\sigma|)=1$, and $|\sigma|/3 \leq \dim_{k} \psi \leq |\sigma|$ for $|\sigma|$ divisible by 3. 
\end{theorem}

\begin{proof} 
Let $M$ be the $m$-dimensional simple left $S$-module corresponding to $\psi$. By \cite[Lemma~4.1]{SmithStaniszkis} and the remarks after that result, $M$ is the quotient of some 1-critical graded $S$-module $N$ with multiplicity less equal to $m$. Here, the {\it multiplicity} of $N$ is defined as $[(1-t)H_N(t)]|_{t=1}$ where $H_N(t)$ is the Hilbert series of $N$. Note that $N$ is an irreducible object in $S$-qgr, which is also $g$-torsionfree. The equivalence of categories, $\underset{g-\text{torsionfree}}{\text{Irred}(S-\text{qgr})} \sim \text{{\tt Irrep}}_{< \infty}(\Lambda_0)$, from \cite[Theorem 7.5]{ATV2} implies that $N$ corresponds the object $N[g^{-1}]_0$ in {\tt Irrep}$_{< \infty}(\Lambda_0)$. Furthermore, $\dim_k N[g^{-1}]_0 = \text{mult}(N)$ as follows. 

Since $N$ is 1-critical, we know by \cite[$\S$2]{ATV2} that the Hilbert series $H_N(t)=f(t)(1-t)^{-1}$ for some $f(t) \in \Z[t,t^{-1}]$. Expanding $H_N(t)$, we see that $\dim_k (N_j) = \text{mult}(N)$ for $j\gg0$. Recall that $g$ is homogeneous of degree 3. Note that $N[g^{-1}]_0 \cdot g^i \subseteq N_{3i}$, and moreover that $\dim_k(N_{3i}) = \text{mult}(N)$ for $i\gg0$. Hence such an $i$:
$$\dim_k N[g^{-1}]_0 ~=~ \dim_k (N[g^{-1}]_0 \cdot g^i) ~\leq~ \dim_k N_{3i} ~=~ \text{mult}(N).$$
Conversely, $N_{3i} \cdot g^{-i}  \subseteq N[g^{-1}]_0$. Hence for $i\gg0$:
$$\text{mult}(N) ~=~ \dim_k(N_{3i}) ~=~ \dim_k(N_{3i} \cdot g^{-i}) ~\leq~ \dim_k N[g^{-1}]_0.$$
Thus, $\dim_k N[g^{-1}]_0 = \text{mult}(N)$ as desired.

Since $\Lambda_0$ is Azumaya, the $\Lambda_0$-module $N[g^{-1}]_0$ has dimension equal to PIdeg($\Lambda_0$) \cite{Artin:azumaya}. Therefore, PIdeg($\Lambda_0$) = mult($N$) $\leq m$. On the other hand, we know by Proposition \ref{prop:dimnsimple} that $m \leq p$. Hence
$$\text{PIdeg}(\Lambda_0) ~\leq~ \dim_k M ~\leq~ p.$$

In the case of $(3, |\sigma|) = 1$, Proposition \ref{prop:pideglam0} implies that PIdeg$\Lambda_0$ = $p$. Thus, $\dim_k M = p$ in this case. For $|\sigma|$ divisible by 3, we also know from Proposition~\ref{prop:pideglam0} that $p/3 \leq \dim_k M \leq p$. We conclude the result by applying Corollary \ref{cor:PIdegS}.
\end{proof}

\begin{corollary} In the case of $|\sigma| < \infty$, a generic finite-dimensional, $g$-torsionfree irrep of $S$ has dimension equal to $|\sigma|$. 
\end{corollary}

\begin{proof}
Since PIdeg($S$) = PIdeg($S[g^{-1}]$), this is a standard consequence of Corollary \ref{cor:PIdegS} and Proposition \ref{prop:dimnsimple}.
\end{proof}

Finally in the case of $|\sigma|<\infty$, we study $g$-torsion finite-dimensional irreducible representations of three-dimensional Sklyanin algebras. In other words, we determine the dimensions of irreps of the twisted homogeneous coordinate rings $B$ arising in Theorem \ref{thm:ATVSkly3}.

\begin{theorem} \label{thm:fin.ord.g-tors} For $n:=|\sigma| < \infty$, a non-trivial $g$-torsion finite-dimensional irrep of $S$ has dimension equal to the PI degree of $B$. 
\end{theorem}

\begin{proof} Recall that  irreps of $A$ bijectively correspond to simple left $A$-modules. For a graded ring $A = \bigoplus_{i \in \N} A_i$, let {\tt Irrep}$_{< \infty}^o A$ denote the set of simple finite-dimensional left $A$-modules that are not annihilated by the irrelevant ideal $A_+ = \bigoplus_{i \geq 1} A_i$. 
Notice that
$\underset{g-\text{torsion}}{\text{{\tt Irrep}}_{< \infty} S}$  = $\text{{\tt Irrep}}_{< \infty} B.$
Thus, it suffices to show that all modules in $\text{{\tt Irrep}}_{< \infty}^o B$ have maximal dimension (which is equal to PIdeg($B$) by Proposition \ref{prop:dimnsimple}). We proceed by establishing the following claims.
\medskip

\begin{claim} \label{claim1} We can reduce the task of studying $\text{{\tt Irrep}}_{< \infty}^o B$ to studying $\text{{\tt Irrep}}_{< \infty}^o \overline{B}$ for $\overline{B}$ some factor of $B$. Furthermore,
$\text{{\tt Irrep}}_{< \infty}^o \overline{B} = \text{{\tt Irrep}}_{< \infty}^o C$
for $C=B(Y_n, \mc{L}|_{Y_n}, \sigma|_{Y_n})$, the twisted homogeneous coordinate ring of an irreducible $\sigma$-orbit of $E$. 
\end{claim}

\begin{claim} \label{claim2} The ring $C$ is isomorphic to a graded matrix ring.
\end{claim}

\begin{claim} \label{claim3} The modules of $\text{{\tt Irrep}}_{< \infty}^o C$ all have maximal dimension, which is equal to PIdeg($C$) = $|\sigma|$.
\end{claim}

\noindent {\it Proof of Claim \ref{claim1}:}
Take $M \in \text{{\tt Irrep}}_{< \infty}^o B$ and let $P$ be the largest graded ideal contained in ann$_B M$, which itself is a prime ideal. Set $\overline{B} = B/P$.
We have by Lemma \ref{lem:GKdimA/P} that GKdim($\overline{B}$)=1, so
Kdim($\overline{B}$) = 1 \cite[8.3.18]{McConnellRobson}. We also have by \cite[Lemma 4.4]{ArtinStafford} that the height one prime $P$ of $B$ corresponds to an $\sigma$-irreducible maximal closed subset of $E$. Namely, since $\sigma$ is given by translation, $P$ is associated to a $\sigma$-orbit of $|\sigma|$ points of $E$, an orbit denoted by  $Y_n$. Here, $n:= |\sigma|$. Now we have a natural homomorphism:
$$\phi: B \lra B(Y_n, \mc{L}|_{Y_n}, \sigma|_{Y_n}) =: C$$
given by restrictions of sections, with ker($\phi$) = $P$. The map $\phi$ is also surjective in high degree, whence $\overline{B}$ is isomorphic to $C$ in high degree.
 
Since $\overline{B}$ and $C$ are PI, we have by Kaplansky's theorem \cite[Theorem~13.3.8]{McConnellRobson} that for each of these rings: the maximal spectrum equals the primitive spectrum. Thus, it suffices to verify the following subclaim.

\begin{subclaim} \label{claim4} For a graded $k$-algebra $A = \bigoplus_{i \geq 0} A_i$, let max$^o A$ denote the set of maximal ideals of $A$ not containing the irrelevant ideal $A_+$. Then, there is a bijective correspondence between
$\max^o C$ and $\max^o \overline{B}$,
given by $I \mapsto \overline{B} \cap I$. Moreover $\overline{B}/(\overline{B} \cap I) \cong C/I$. 
\end{subclaim}

\noindent {\it Proof of Subclaim \ref{claim4}:} We know that $\overline{B}_{\geq m} = C_{\geq m}$ for $m \gg0$. Since $B$ is generated in degree 1, we have for any such $m$:
$$J:= (\overline{B}_+)^m = \overline{B}_{\geq m} = C_{\geq m} \supseteq (C_+)^m,$$
is a common ideal of $\overline{B}$ and $C$. (The last inclusion may be strict as $C$ need not be generated in degree 1.) 

For one inclusion, let $I \in \max^o C$. We want to show that $\overline{B} \cap I \in \max^o \overline{B}$. Note that $\overline{B}/(\overline{B} \cap I) \cong (\overline{B} + I)/I$ as rings. Furthermore,
$$C ~\supseteq~ \overline{B}+I ~\supseteq~ J+I ~\supseteq~ (C_+)^m + I = C.$$
The last equality is due to $C_+$ and $I$ being comaximal in $C$. Hence $C=\overline{B}+I$ and $\overline{B}/(\overline{B} \cap I) \cong C/I$ is a simple ring. Consequently, $\overline{B} \cap I \in \max^o \overline{B}$.

Conversely, take $M \in \max^o \overline{B}$ and we want to show that $M = \overline{B} \cap Q$ for some $Q \in \max^o C$. We show that $CMC \neq C$. Suppose not. Then, 
$$(\overline{B}_+)^{2m} = J^2 = JCJ = JCMCJ = JMJ \subseteq M,$$
which implies that $\overline{B}_+ \subseteq M$ as $M$ is prime. This is a contradiction. Now $CMC$ is contained in some maximal ideal $Q$ of $C$. Since $M \in \max^o \overline{B}$, we have that $Q \in \max^o C$; else $M \subseteq Q \cap \overline{B} = \overline{B}_+$. By the last paragraph, we know that $\overline{B} \cap Q \in \max^o \overline{B}$ with
$\overline{B} \cap Q = M \in \max \overline{B}$. 

Thus, Claim \ref{claim1} is verified and so it suffices to examine $\text{{\tt Irrep}}_{< \infty}^o C$. 
\medskip

\noindent {\it Proof of Claim \ref{claim2}:} Recall that $n := |\sigma|$. 
We will show that $C$ is isomorphic to the graded matrix ring:
 
{\footnotesize
$$R := \left( \begin{array}{ccccc}
              T       & x^{n-1} T & x^{n-2}T   & \dots  & xT\\
              xT      & T         & x^{n-1}T   & \dots  & x^2 T\\
              x^2 T   & xT        & T          & \dots  & x^3 T\\
              \vdots  & \vdots    & \vdots     & \ddots & \vdots\\
            x^{n-1} T & x^{n-2} T & x^{n-3}T   & \dots  & T
              \end{array} \right)$$}
              
\noindent with $T = k[x^n]$. Say $Y_n := \{p_1, \dots, p_n\}$, the $\sigma$-orbit of $n$ distinct points $p_i$. Let $\mc{L}' := \mc{L}|_{Y_n}$ and $\sigma' := \sigma|_{Y_n}$. Reorder the $\{p_i\}$ to assume that $\sigma'(p_i) = p_{i+1}$ and $\sigma'(p_n) = p_1$ for $1 \leq i \leq n-1$. By Definition \ref{def:thcr}, $C = \bigoplus_{d \geq 0} H^0(Y_n, \mc{L}'_d)$ where $\mc{L}'_0 = \bigoplus_{i=1}^n \mc{O}_{p_i}$, and $\mc{L}'_1 = \bigoplus_{i=1}^n \mc{O}_{p_i}(1)$, and 
$$\mc{L}'_d = \mc{L}' \ten_{\mc{O}_{Y_n}} (\mc{L}')^{\sigma'} \ten_{\mc{O}_{Y_n}} \dots \ten_{\mc{O}_{Y_n}} (\mc{L}')^{{\sigma'}^{d-1}}.$$
Note that $$(\mc{L}')^{\sigma} = \mc{O}_{p_n}(1) \oplus \mc{O}_{p_1}(1) \oplus \dots \oplus \mc{O}_{p_{n-1}}(1) \cong \mc{L}'.$$ 
Therefore, $\mc{L}'_d \cong (\mc{L}')^{\ten d}$. If $i=j$, then $\mc{O}_{p_i}(1) \ten \mc{O}_{p_j}(1) \cong \mc{O}_{p_i}(2)$. Otherwise for $i \neq j$, the sheaf $\mc{O}_{p_i}(1) \ten \mc{O}_{p_j}(1)$ has empty support. 
Hence $(\mc{L}')^{\ten d} \cong \bigoplus_{i=1}^n \mc{O}_{p_i}(d)$. So, the ring $C$ shares the same $k$-vector space structure as a sum of $n$ polynomial rings in one variable, say $k[u_1] \oplus \dots \oplus k[u_n]$.

Next, we study the multiplication of the ring $C$. Write $C_i$ as $k u_1^i \oplus \dots \oplus k u_n^i$. The multiplication $C_i \times C_j \ra C_{i+j}$ is defined on basis elements $(u_1^k, \dots, u_n^k)$ for $k=i,j$, which will extend to $C$ by linearity. All of the following sums in indices are taken modulo $n$. 

Observe that
\[
\begin{array}{rl}
C_i \times C_j &= H^0(Y_n, \mc{L}'_i) \ten H^0(Y_n, \mc{L}'_j)\\
               &= H^0(Y_n, \mc{L}'_i) \ten H^0(Y_n, (\mc{L}'_j)^{(\sigma')^i}).
\end{array}
\] 
The multiplicative structure presented in Definition \ref{def:thcr} implies that the multiplication of $C$ is given as 
$$(u_1^i, \dots, u_n^i) \ast (u_1^j, \dots, u_n^j)
= (u_1^i u_{1-i}^j, \dots, u_n^i u_{n-i}^j).$$

We now show that $C$ is isomorphic to the graded matrix ring $R$. Define a map $\phi: C \ra R$ by 
$$(u_1^i, \dots, u_n^i) ~\mapsto~ \left( \text{Row}(u_l^i) \right)_{l=1,\dots,n} =: {\tt M_i} ~\in R.$$
                                       
\noindent Here Row($u_l^i$) denotes the row $l$ of a matrix with the entry $u_l^i$ in column $l-i$ modulo $n$, and zeros entries elsewhere. In other words, the matrix ${\tt M_i}$ has a degree $i$ entry in positions $(l, l-i)$ for $l=1,\dots,n$ and zeros elsewhere. We see that ${\tt M_i} \in R$.

The map $\phi$ is the ring homomorphism as follows. First note that
\[
\begin{array}{rl}
\phi((u_1^i, \dots, u_n^i) \ast (u_1^j, \dots, u_n^j)) &= \phi((u_1^i u_{1-i}^j, \dots, u_n^i u_{n-i}^j))\\
&= \left( \text{Row}(u_l^i u_{l-i}^j) \right)_{l=1,\dots,n}.
\end{array}
\]
Here, the entry $u_l^i u_{l-i}^j$ appears in positions $(l, l-(i+j))$ for $l=1,\dots,n$.

On the other hand, 
\[
\phi((u_1^i, \dots, u_n^i)) \cdot \phi((u_1^j, \dots, u_n^j)) ~=~ {\tt M_i} \cdot {\tt M_j}.
\]
The nonzero entries of ${\tt M_i}$, namely $\{u_l^i\}$, appear in positions $(l, l-i)$ for $l=1,\dots,n$. Whereas the nonzero entries of ${\tt M_j}$, the $\{u_l^j\}$, are in positions $(l, l-j) = (l-i, l-(i+j))$ for $l,i=1,\dots,n$. Therefore, the product ${\tt M_i} \cdot {\tt M_j}$ has nonzero entries in positions $(l, l-(i+j))$; these entries are equal to $u_l^i u_{l-i}^j$ for $l = 1,\dots,n$.
Thus, we have a ring homomorphism $\phi$ between $C$ and $R$ which is bijective. This concludes the proof of Claim \ref{claim2}. 
\medskip

\noindent {\it Proof of Claim \ref{claim3}:}
We aim to show that $M \in$ {\tt Irrep}$_{<\infty}^o R$ has maximal $k$-vector space dimension, which is equal to $n$. Let $s$ denote the diagonal matrix, diag($x^n$) in $R$, and note that $sR = R_+$. Consider the ring $R[s^{-1}]$. For $M \in$ {\tt Irrep}$_{<\infty}^o R$, construct the module $$0 \neq M[s^{-1}] = R[s^{-1}] \ten_R M ~~\in R[s^{-1}]-\text{mod}.$$ Since $M$ is simple, $M[s^{-1}]$ is also simple by \cite[Corollary 10.16]{GoodearlWarfield}. Furthermore, $\dim_k M[s^{-1}]$ equals $\dim_k M$ as $M$ is $s$-torsionfree \cite[Exercise 10K]{GoodearlWarfield}. Thus, we have an inclusion of sets {\tt Irrep}$_{< \infty}^o R \subseteq$ {\tt Irrep}$_{<\infty} R[s^{-1}]$ given by $M \mapsto M[s^{-1}]$. 
Now observe that:

{\footnotesize
$$R[s^{-1}] = \left( \begin{array}{ccccc}
                     T' & x^{n-1} T' & & xT'\\
                     xT' & T'    &  & x^2 T'\\
                     \vdots & \vdots & \ddots & \vdots\\
                     x^{n-1} T' & x^{n-2} T' &  & T'
                     \end{array} \right)$$}
                     
\noindent for $T' = k[(x^n)^{\pm 1}]$. So $R[s^{-1}] \subseteq$ Mat$_n(T')$ as rings, i.e. the rings have the same multiplicative structure. The rings are also equal as sets. Hence $R[s^{-1}]$ =  Mat$_n(T')$ as rings. Since all simple modules of Mat$_n(T')$ have dimension $n$, we have verified Claim \ref{claim3}.

This concludes the proof of the theorem. 
\end{proof}

\begin{corollary} \label{cor:PIdegB}
Given a three-dimensional Sklyanin algebra $S(a,b,c)$ with $|\sigma_{abc}|<\infty$, and $S/Sg \cong$ a twisted homogeneous coordinate ring $B$, we have that PIdeg($S$) = PIdeg($B$) = $|\sigma|$.
\end{corollary}

\begin{proof}
In the proof of the preceding proposition, $|\sigma|$ is equal to the maximal dimension of the finite-dimensional irreps of the rings $R$ and $C$. By Claim \ref{claim1},  finite-dimensional irreps of $R$ and $C$ correspond to modules in {\tt Irrep}$_{< \infty}^o B$ of the same dimension, so $|\sigma|$ = PIdeg($B$) (Proposition \ref{prop:dimnsimple}). Moreover, PIdeg($S$) = $|\sigma|$ by Corollary \ref{cor:PIdegS}.
\end{proof}


\section{Irreducible representations of $Skly_3$ of intermediate dimension}

In this section, we consider three-dimensional Sklyanin algebras equipped with automorphism $\sigma$ so that $|\sigma| < \infty$ is divisible by 3. We show that the lower bound of dimensions of $g$-torsionfree irreps of $Skly_3$, namely $|\sigma|/3$, can be achieved; refer to Theorem \ref{thm:fin.ord.g-tfree}. Here, we employ the Sklyanin algebra $S(1,-1,-1)$, which we denote by $\hat{S}$.

\begin{prop} \label{prop:Simp2} The automorphism $\sigma_{1,-1,-1}$, corresponding to the Sklyanin algebra $\hat{S}=S(1,-1,-1)$,
 has order 6. Hence, the PI degree of $\hat{S}$ is also 6. Moreover, there exist 2-dimensional irreducible representations of $\hat{S}$.
\end{prop}

\begin{proof} 
For the proof of the first two statements, apply Proposition \ref{prop:sigmafinite} and Corollary~\ref{cor:PIdegS}. Now, we construct 2-dimensional irreps of $\hat{S}$ with the assistance of the computer algebra software Maple. If $\pi \in$ {\tt Irrep}$_2 \hat{S}$, then for the generators $x$, $y$, $z$ of $\hat{S}$ let $\pi(x)=:X$, $\pi(y)=:Y$, and $\pi(z)=:Z$. 

Without loss of generality, we assume that $X$ is in Jordan-Canonical form.
If $X$ is non-diagonal, then $Y$ and $Z$ are also upper triangular, and this yields reducible representations.
Hence, we assume that $X$ is diagonal. One gets five sets of representations, three of which are reducible. The irreducible representations are given as follows:

{\footnotesize
\[
X = \left(\begin{array}{cc}
         (\pm i +1)z_4 & 0\\
         0              & -(\pm i -1)z_4
          \end{array} \right),          
~Y = \left(\begin{array}{cc}
         z_4 & \frac{\pm i z_4^2}{z_3}\\
          \pm i z_3    & z_4
          \end{array} \right),
~Z = \left(\begin{array}{cc}
         z_4 & -\frac{z_4^2}{z_3}\\
          z_3    & z_4
          \end{array} \right)
\]
}
\noindent and
{\footnotesize
\[
X = \left(\begin{array}{cc}
         \xi(\xi^2 - \xi -1)z_4 & 0\\
         0              & \frac{(\xi^2 - \xi + 1)z_4}{\xi -1}
          \end{array} \right),          
~Y = \left(\begin{array}{cc}
         (\xi^2 - 1)z_4 & \frac{\xi z_4^2}{z_3}\\
          \xi z_3    & (\xi^2 -1)z_4
          \end{array} \right),
~Z = \left(\begin{array}{cc}
         z_4 &- \frac{z_4^2}{z_3}\\
          z_3    & z_4
          \end{array} \right)
\]
}

\noindent where $\xi$ is a primitive twelve root of unity, and $z_3, z_4$ are arbitrary scalars in $k$.
\end{proof}

\begin{remark}
Recall from $\S$2.2 that nontrivial irreps of PI rings have nontrivial central annihilators by Schur's Lemma, and hence the center of $\hat{S}$ controls the behavior of {\tt Irrep}$_{< \infty} \hat{S}$. Note that $Z(\hat{S})$ is generated by three elements of degree 6 and one of degree 3 by Proposition~\ref{prop:Simp2} and \cite[Theorem~3.7]{SmithTate}. We use Affine, a noncommutative algebra package of Maxima, to compute these generators:
\[
\begin{array}{rl}
g    &= x^3 - yxz \text{\hspace{.4in}(from Equation \ref{eq:g})},\\
u_1  &= yx^2yxz + xyxyxz + x^2yxyz + x^3yx^2 + x^4yx,\\
u_2  &= yx^4y - xyx^2yx + x^2yxyx - 2x^2yx^2y - x^3yxy,\\
u_3  &= x^3yxz.
\end{array}
\]
There is a degree 18 relation amongst the algebraically independent elements $\{u_i\}$ and the element $g$; this follows from \cite[Theorem 4.7]{SmithTate}.
\end{remark}


\section{Further results}
\label{sec:further}

\subsection{Sklyanin algebras that are module finite over their center}

Three-dimensional Sklyanin algebras that are module finite over their center (or PI Sklyanin algebras) are of particular interest in this work. Recall that $S(a,b,c)$ is PI if and only if $|\sigma_{abc}|<\infty$ \cite[Theorem 7.1]{ATV2}. We know that the PI degree of such $Skly_3$ = $S(a,b,c)$ is, in fact, equal to $|\sigma_{abc}|$ (Corollary \ref{cor:PIdegS}). On the other hand,  we proved in $\S$3 that we have nontrivial matrix solutions to Equations \ref{eq:S(abc)}  precisely when $|\sigma_{abc}|<\infty$, most of which are of dimension $|\sigma| \times |\sigma|$.  Naturally, one would like to classify parameters $a$, $b$, $c$ so that $\sigma_{abc}$ has finite order; this is equivalent to describing the PI Sklyanin algebras, or those $S(a,b,c)$ with nontrivial representation theory.

We do not necessarily assume that we have condition (2) of Definition \ref{def:Skl3}, so $E_{abc}$ could be $\mathbb{P}^2$ or a triangle, for instance (see Equation \ref{eq:E(abc)}). Moreover, by \cite[$\S$1]{ATV1}, $\sigma_{abc}$ is an automorphism of $E_{abc}$ induced by the shift functor on point modules of $S$, which can be explicitly given as follows: 
\begin{equation} \label{eqn:sigma}
\sigma([x:y:z]) = [acy^2 - b^2xz ~:~ bcx^2  - a^2yz ~:~ abz^2 - c^2xy].
\end{equation}

\begin{prop} \label{prop:sigmafinite} Let $\zeta$ be a third root of unity. Then, we have the following results pertaining to $\sigma_{abc}$ of finite order.
\smallskip

{\small
$\bullet$ $|\sigma|=1 \iff$ \text{$[a:b:c] = [1:-1:0]$}, \text{~the origin of~} $E$.
}

{\small
$\bullet$ $|\sigma|=2 \iff$ \text{$[a:b:c] = [1:1:c]$}, \text{~for~} $c \neq 0$, ~$c^3 \neq 1$.
}
\vspace{.2in}

{\small
$\bullet$ $|\sigma|=3 \iff$ \text{$[a:b:c]$} = $\begin{cases}
                                              [1:0:\omega], \text{~for~} \omega = -1, e^{\pi i/3}, e^{5 \pi i/3};\\
                                              [1:\omega:0], \text{~for~} \omega = e^{\pi i/3}, e^{5 \pi i/3}.
                                              \end{cases}$
}

\vspace{.2in}

{\small
$\bullet$ $|\sigma|=4 \iff$ \text{$[a:b:c]$} = $\begin{cases}
                                              \left[ 1:b:\left(\frac{b(b^2+1)}{b+1}\right)^{1/3} \zeta \right], \text{~for~} b \neq 0, -1, \pm i;  b^3 \neq 1;\\
                                              \left[ 1: \left( \frac{f(c)^{2/3}-12(1-c^3)}{6 f(c)^{1/3}} \right) \zeta :c \right], \text{for~} c \neq 0, c^3 \neq 1.
                                              \end{cases}$

\vspace{.1in}
                                              
\indent \indent where ~$f(c) = 108c^3 + 12(12-36c^3+117c^6-12c^9)^{1/3}$.
}

\vspace{.2in}

{\small
$\bullet$ $|\sigma|=5 \iff$ \text{$[a:b:c]$} = $\begin{cases}
                                              \left[ 1:b:\left(\frac{r \pm s^{1/2}}{2b}\right)^{1/3} \zeta \right],   \text{~for~} b \neq 0, ~ b^3 \neq 1; \\
                           \text{\hspace{0.65in}}          b \neq \text{~primitive~} 10^{th} \text{~root of unity},\\
                                              \left[ 1: g(c) :c \right],  \text{~for~} c \neq 0, c^3 \neq 1.
                                              \end{cases}$\\{}\\
\indent \indent where ~$r = -b^5+b^4+b^3+b^2+b-1$ and $s=(b^2-3b+1)(b-1)^2(b^2+b-1)^3$.

\indent Moreover ~$g(c)$~ is ~a ~root ~of

\indent \indent  $Z^6 +(c^3-1)Z^5 + (1-c^3)Z^4 + (-1-c^3)Z^3 + (1-c^3)Z^2 + (c^6-c^3)Z +c^3$.
}

\vspace{.2in}

{\small
$\bullet$ $|\sigma|=6 \iff$ \text{$[a:b:c]$} = $\begin{cases}
                                              \left[ 1:b:\zeta \right], \left[ 1:b:b \zeta \right], \text{~for~} b \neq 0, b^3 \neq 1,\\
                                              \left[ 1: \zeta :c \right], \left[ 1: c \zeta :c \right], \text{~for~} c \neq 0, c^3 \neq 1, \zeta \neq 1.
                                              \end{cases}$
}
\end{prop}

\smallskip

\begin{proof}
This is achieved via a Maple computation; see \cite[Proposition A.1]{Walton:thesis} for details.
\end{proof}

We point out that a conjecture of Artin-Schelter on the PI degree of a PI Sklyanin algebra follows from the proposition above.

\begin{corollary} \label{cor:conjAS} Conjecture 10.37(ii) of \cite{ArtinSchelter} holds for $r=3$ (for type A quadratic AS-regular algebras of dimension 3). To say, if $a=b$, then $S$ has rank $2^2$ over its center. Moreover, all higher $n^2$ occur as ranks.
\end{corollary}

\begin{proof}
Refer to the $|\sigma|=2$ case of Proposition \ref{prop:sigmafinite}, and apply Corollary \ref{cor:PIdegS}, respectively.
\end{proof}

\begin{remark}
In \cite[$\S$4]{Berenstein}, it is stated that the algebra $S(1,1,c)$ is of dimension 9 over its center. However by Corollary \ref{cor:PIdegS} and the previous proposition, the algebra $S(1,1,c)$ has PI degree 2. \end{remark}

Observe that for $n=2, 4,5,6$, there are 1-parameter families of solutions for which $|\sigma_{abc}|=n$ and we draw our attention to degeneration limits of some of these families. 
For instance, consider the degeneration limit $c=0$ for the $|\sigma|=2$ case. This yields the algebra
$$S(1,1,0) = k\{x,y,z\}/(yz+zy, ~zx+xz, ~xy+yx),$$
a skew polynomial ring \cite[Chapter 2]{GoodearlWarfield}. 
On the other hand, consider the degeneration limit $b = \pm i$ for the $|\sigma|=4$, which yields the algebra
$$S(1,\pm i,0) = k\{x,y,z\}/(yz \pm izy, ~zx \pm ixz, ~xy \pm iyx),$$
also a skew polynomial ring. 
Both of these degenerate cases arise in results on orbifolds with discrete torsion \cite[$\S$4.1]{BJL}. 

Although the analysis of the geometry of $S(a,b,c)$ with $|\sigma|<\infty$ is the subject of future work, we describe the centers of $S(a,b,c)$ at $(a,b,c) = (1,1,0), (1, \pm i, 0)$ and other similar degeneration limits for now. For the remainder of this subsection, we restrict ourselves to the case where the field $k = \C$.

\begin{definition} \label{def:CY} \cite[Definition 1.4.1]{CoxKatz} \cite[$\S$1]{He:lectures} A {\it Calabi-Yau threefold (CY)} is defined to be a normal quasi-projective algebraic variety $X$ with:
\begin{enumerate}
\item trivial dualizing sheaf $\omega_X \cong \mc{O}_X$;
\item $H^i(X,\mc{O}_X) = 0$ for $i = 1,2$; and 
\item at worst, canonical Gorenstein singularities (see \cite{Reid:young} for more details).
\end{enumerate}
\end{definition}

First, we state a general result about the CY geometry and center of the Sklyanin algebra $S(1,q,0)$, with $q$ a root of unity.

\begin{prop} \label{prop:centerdeg}
Let $q$ be a primitive $n^{\text{th}}$ root of unity. Then, the center $Z$ of $S(1,q,0)$ is generated by $u_1 = x^n$, $u_2=y^n$, $u_3=z^n$ and $g=xyz$, subject to relation $(-g)^n + u_1 u_2 u_3$. Here, $g$ arises as the central cubic element  in Equation \ref{eq:g}. Furthermore, Spec($Z$) is the affine toric CY three-fold: $\C^3/(\Z_n \times \Z_n)$.
\end{prop}

\begin{proof}
The center $Z$ of $S(a,b,c)$ with $|\sigma_{abc}|=n<\infty$ has three algebraically independent generators $u_i$ of degree $n$, and one generator $g$ of degree 3, subject to a relation of degree $3n$ \cite[Theorems 3.7, 4.7]{SmithTate}. 
This applies to the algebra $S=S(1,q,0)$ as the automorphism $\sigma_{1,q,0}$ from Equation \ref{eq:g}
has order $n$. 

Let $(a,b,c) = (1,q,0)$. We have by Equation \ref{eq:g} that $g= (q^3 - 1)xyz$ is a central element of $S(1,q,0)$. Without loss of generality, we can take $g = xyz$.
Furthermore, we verify by brute force that the algebraically independent elements  \{$u_1=x^n$, $u_2=y^n$, $u_3=z^n$\} are also generators of the center of $S(1,q,0)$, and that the relation amongst the generators is given by $(-g)^n + u_1 u_2 u_3$. 

Now, we will show that  $X$:=Spec($Z(S(1,q,0))$) is isomorphic to $\C^3/(\Z_n \times \Z_n)$ as varieties in $\C^4$. Let $G := \Z_n \times \Z_n$, and let $\nu$ be a primitive $n^{th}$ root of unity. Take $G$ to be generated by the matrices, diag($\nu$, 1, $\nu^{-1}$) and diag($\nu$, $\nu^{-1}$, 1), so that $G$ acts on the vector space $\C x \oplus \C y \oplus \C z$. With the aid of the computational algebra system MAGMA, we have that $x^n$, $y^n$, $z^n$, and $g$ are precisely the generators of  the invariant ring $\C[x,y,z]^G$. Hence, $Z \cong \C[x,y,z]^G$, and
$$X ~=~ \text{Spec}(\C[x,y,z]^G) ~=~ \C^3/G ~=~ \C^3/(\Z_n \times \Z_n),$$
an affine  GIT quotient variety \cite[Definition 10.5]{MillerSturmfels}.

Here, we show that $X$ is an affine toric CY threefold, even though this fact is well-known. Condition (2) of Definition \ref{def:CY} holds as $X$ is an affine variety. Moreover, since $G$ is a finite subgroup of $SL(3,\C)$, the dualizing sheaf of $\C^3/(\Z_n \times \Z_n)$ is trivial \cite[$\S$3.1]{Cox:Update}. Thus, condition (1) also holds.
Since $X$ is a hypersurface in affine space, $X$ is Gorenstein. Thus, $X$ has canonical singularities by \cite[Remark 1.8]{Reid:canonical}, and the last condition of Definition \ref{def:CY} holds for $X$. The variety $X \cong \C^3/(\Z_n \times \Z_n)$ is also a toric variety since $\Z_n \times \Z_n$ is an abelian subgroup of $SL(3, \C)$ \cite[$\S$3.3]{Cox:Update}. 
\end{proof}

\begin{corollary} \label{cor:CY} Consider the degeneration limits $c=0$ for the $|\sigma|=2$ case, and $b=\pm i$ for the $|\sigma|=4$ case of Proposition \ref{prop:sigmafinite}. The centers of both $S(1,1,0)$ and $S(1,\pm i, 0)$ produce affine toric CY three-folds:  Spec($Z$) = $\C^3/(\Z_2 \times \Z_2)$ and $\C^3/(\Z_4 \times \Z_4)$, respectively.
\qed
\end{corollary}

\begin{remark} \label{rmk:crepant} A smooth crepant resolution of $\C^3/(\Z_n \times \Z_n)$ is provided in \cite[Theorem~0.1]{Nakamura}, or more explicitly in \cite{CrawReid}.
\end{remark}

\subsection{Classifying irreducible representations of deformed Sklyanin algebras}

As discussed in the introduction, this work is motivated by \cite{BJL} in which the authors study various deformations of the N=4 SYM theory in four dimensions. Such deformed theories include {\it relevant} deformations, equipped with mass or linear terms. Here, the superpotential is $\Phi = \Phi_{marg} + \Phi_{rel}$ given by:
\[
\begin{array}{rl}
\Phi_{marg} &= axyz + byxz + \frac{c}{3}(x^3 + y^3 + z^3), \\
\Phi_{rel}  &= \frac{d_1}{2}x^2 + \frac{d_2}{2}(y^2 +z^2) + (e_1 x + e_2 y + e_3 z).
\end{array}
\]
Thus, the finite-dimensional irreps of the corresponding superpotential algebras are desired. In fact, these algebras denoted by $S_{def}$, are realized as ungraded deformations of $Skly_3$. 

To employ techniques from noncommutative projective algebraic geometry, we observe that the $S_{def}$ are homomorphic images of graded {\it central extensions} $D$ of $Skly_3$ in the sense of \cite{LSV}.
Furthermore, finite-dimensional irreps of $D$ are precisely those of $Skly_3$ or of $S_{def}$.
Hence, the analysis of {\tt Irrep}$_{< \infty} S_{def}$ falls into two tasks: the study of the representation theory of each of the graded algebras $Skly_3$ and $D$. As we have completed the first task, we leave the study of the set {\tt Irrep}$_{< \infty} D$ to the reader. Partial results are available in \cite[Chapter 5]{Walton:thesis}.

\subsection{Classifying irreducible representations of degenerate Sklyanin algebras}

In this section, we consider twelve algebras that were omitted from the family of algebras $S(a,b,c)$ in Definition \ref{def:Skl3}, so-called {\it degenerate Sklyanin algebras}.

\begin{definition} \label{def:Sdeg} The rings $S(a,b,c)$ from Definition \ref{def:Skl3} with $[a:b:c] \in \mf{D}$ are called {\it degenerate Sklyanin algebras}, denoted by $S_{deg}$\index{S_{deg}}.
\end{definition}

It is shown in \cite{Walton:Sdeg} that $S_{deg}$ has an entirely different ring-theoretic and homological structure than its counterpart $Skly_3$. Likewise, we will see in the proposition below that the representation theory of $S_{deg}$ does not resemble that of $Skly_3$ as we construct irreps of $S_{deg}$ for every even dimension. 

\begin{prop} Each degenerate Sklyanin algebra has irreducible representations of dimension $2n$ for all positive integers $n$. 
\end{prop}

\begin{proof}
Let $I_n$ be the $n \times n$ identity matrix, let $J_n$ be the $n \times n$ anti-identity matrix, and let {\bf 0} be the $n \times n$ zero matrix.
The following triples of matrices ($X$, $Y$, $Z$) are irreps of the algebras $S(1,0,0)$, $S(0,1,0)$, and $S(0,0,1)$, respectively:

{\footnotesize
\[
X=\left(\begin{array}{cc}
        x\cdot I_n & {\bf 0}\\
        {\bf 0}& {\bf 0}
        \end{array}\right),~~
Y=\left(\begin{array}{cc}
        {\bf 0}&{\bf 0}\\
        {\bf 0} & y\cdot I_n
        \end{array}\right),~~
Z=\left(\begin{array}{cc}
       {\bf 0} & z\cdot J_n\\
        {\bf 0} & {\bf 0}
        \end{array}\right);
\]
}

\vspace{-.2in}

{\footnotesize
\[       
X=\left(\begin{array}{cc}
        x\cdot I_n & {\bf 0}\\
        {\bf 0}& {\bf 0}
        \end{array}\right),~~
Y=\left(\begin{array}{cc}
        {\bf 0} & y\cdot J_n\\
        {\bf 0} & {\bf 0}
        \end{array}\right),~~
Z=\left(\begin{array}{cc}
       {\bf 0}&{\bf 0}\\
        {\bf 0} & z\cdot I_n
        \end{array}\right);
\]
}

\vspace{-.2in}

{\footnotesize
\[
X=\left(\begin{array}{cc}
        {\bf 0}& x\cdot J_n\\
        {\bf 0}& {\bf 0}
        \end{array}\right),~~
Y=\left(\begin{array}{cc}
        {\bf 0}& {\bf 0}\\
        y\cdot J_n & {\bf 0}
        \end{array}\right),~~
Z=\left(\begin{array}{cc}
       {\bf 0} & z\cdot J_n\\
        {\bf 0} & {\bf 0}
        \end{array}\right).
\]
}

On the other hand, we have the following irreps for the degenerate Sklyanin algebras $S(1,b,c)$ with $b^3 = c^3 =1$; here $p$ is arbitrary in $k^{\times}$:

{\footnotesize
\[
X=\left(\begin{array}{cc}
        -bx\cdot I_n & {\bf 0}\\
        {\bf 0}& x\cdot I_n
        \end{array}\right),~~
Y=\left(\begin{array}{cc}
        {\bf 0}& -\frac{b^2c y}{p}\cdot J_n\\
        {\bf 0}& {\bf 0}
        \end{array}\right),~~
Z=\left(\begin{array}{cc}
       {\bf 0} & {\bf 0}\\
        \frac{px^2}{y} \cdot J_n & {\bf 0}
        \end{array}\right).
\]
}

\vspace{-.35in}
\end{proof}

\bibliography{biblio}

\def\cprime{$'$}\def\cprime{$'$}\def\cprime{$'$}\def\cprime{$'$}\def\cprime{$'$}\def\cprime{$'$}\def\cprime{$'$}
\begin{thebibliography}{10}

\bibitem{Artin:azumaya}
M.~Artin.
\newblock On {A}zumaya algebras and finite dimensional representations of
  rings.
\newblock {\em J. Algebra}, 11:532--563, 1969.

\bibitem{Artin}
M.~Artin.
\newblock Geometry of quantum planes.
\newblock In {\em Azumaya algebras, actions, and modules ({B}loomington, {IN},
  1990)}, volume 124 of {\em Contemp. Math.}, pages 1--15. Amer. Math. Soc.,
  Providence, RI, 1992.

\bibitem{ArtinDeJong}
M.~Artin and A.~J. de~Jong.
\newblock Stable orders over surfaces (course notes).
\newblock Math 711, Winter 2004, University of Michigan.

\bibitem{ArtinSchelter}
M.~Artin and W.~F. Schelter.
\newblock Graded algebras of global dimension {$3$}.
\newblock {\em Adv. in Math.}, 66(2):171--216, 1987.

\bibitem{ArtinStafford}
M.~Artin and J.~T. Stafford.
\newblock Noncommutative graded domains with quadratic growth.
\newblock {\em Invent. Math.}, 122(2):231--276, 1995.

\bibitem{ATV1}
M.~Artin, J.~Tate, and M.~Van~den Bergh.
\newblock Some algebras associated to automorphisms of elliptic curves.
\newblock In {\em The {G}rothendieck {F}estschrift, {V}ol.\ {I}}, volume~86 of
  {\em Progr. Math.}, pages 33--85. Birkh\"auser Boston, Boston, MA, 1990.

\bibitem{ATV2}
M.~Artin, J.~Tate, and M.~Van~den Bergh.
\newblock Modules over regular algebras of dimension {$3$}.
\newblock {\em Invent. Math.}, 106(2):335--388, 1991.

\bibitem{Berenstein}
D.~Berenstein.
\newblock Reverse geometric engineering of singularities.
\newblock {\em J. High Energy Phys.}, (4):No. 52, 18, 2002.

\bibitem{BJL}
D.~Berenstein, V.~Jejjala, and R.~G. Leigh.
\newblock Marginal and relevant deformations of {$N=4$} field theories and
  non-commutative moduli spaces of vacua.
\newblock {\em Nuclear Phys. B}, 589(1-2):196--248, 2000.

\bibitem{Bocklandt:CY3}
R.~Bocklandt.
\newblock Graded {C}alabi {Y}au algebras of dimension 3.
\newblock {\em J. Pure Appl. Algebra}, 212(1):14--32, 2008.

\bibitem{BLS}
R.~Bocklandt, L.~Le~Bruyn, and S.~Symens.
\newblock Isolated singularities, smooth orders, and {A}uslander regularity.
\newblock {\em Comm. Algebra}, 31(12):6019--6036, 2003.

\bibitem{BrownGoodearl}
K.~A. Brown and K.~R. Goodearl.
\newblock Homological aspects of {N}oetherian {PI} {H}opf algebras of
  irreducible modules and maximal dimension.
\newblock {\em J. Algebra}, 198(1):240--265, 1997.

\bibitem{book:BrownGoodearl}
K.~A. Brown and K.~R. Goodearl.
\newblock {\em Lectures on algebraic quantum groups}.
\newblock Advanced Courses in Mathematics. CRM Barcelona. Birkh\"auser Verlag,
  Basel, 2002.

\bibitem{Cox:Update}
D.~A. Cox.
\newblock Update on toric geometry.
\newblock In {\em Geometry of toric varieties}, volume~6 of {\em S\'emin.
  Congr.}, pages 1--41. Soc. Math. France, Paris, 2002.

\bibitem{CoxKatz}
D.~A. Cox and S.~Katz.
\newblock {\em Mirror symmetry and algebraic geometry}, volume~68 of {\em
  Mathematical Surveys and Monographs}.
\newblock American Mathematical Society, Providence, RI, 1999.

\bibitem{CrawReid}
A.~Craw and M.~Reid.
\newblock How to calculate {$A$}-{H}ilb {$\Bbb C^3$}.
\newblock In {\em Geometry of toric varieties}, volume~6 of {\em S\'emin.
  Congr.}, pages 129--154. Soc. Math. France, Paris, 2002.

\bibitem{GoodearlWarfield}
K.~R. Goodearl and R.~B. Warfield, Jr.
\newblock {\em An introduction to noncommutative {N}oetherian rings}, volume~61
  of {\em London Mathematical Society Student Texts}.
\newblock Cambridge University Press, Cambridge, second edition, 2004.

\bibitem{He:lectures}
Y.~H. He.
\newblock Lectures on {D}-branes, gauge theories and {C}alabi-{Y}au
  singularities (preprint). http://arxiv.org/abs/hep-th/0408142.

\bibitem{KrauseLenagan}
G.~R. Krause and T.~H. Lenagan.
\newblock {\em Growth of algebras and {G}elfand-{K}irillov dimension},
  volume~22 of {\em Graduate Studies in Mathematics}.
\newblock American Mathematical Society, Providence, RI, revised edition, 2000.

\bibitem{LeBruyn:centsing}
L.~Le~Bruyn.
\newblock Central singularities of quantum spaces.
\newblock {\em J. Algebra}, 177(1):142--153, 1995.

\bibitem{LSV}
L.~Le~Bruyn, S.~P. Smith, and M.~Van~den Bergh.
\newblock Central extensions of three-dimensional {A}rtin-{S}chelter regular
  algebras.
\newblock {\em Math. Z.}, 222(2):171--212, 1996.

\bibitem{LeighStrassler}
R.~G. Leigh and M.~J. Strassler.
\newblock Exactly marginal operators and duality in four-dimensional {$N=1$}
  supersymmetric gauge theory.
\newblock {\em Nuclear Phys. B}, 447(1):95--133, 1995.

\bibitem{Maldacena}
J.~Maldacena.
\newblock The large {$N$} limit of superconformal field theories and
  supergravity.
\newblock {\em Adv. Theor. Math. Phys.}, 2(2):231--252, 1998.

\bibitem{McConnellRobson}
J.~C. McConnell and J.~C. Robson.
\newblock {\em Noncommutative {N}oetherian rings}, volume~30 of {\em Graduate
  Studies in Mathematics}.
\newblock American Mathematical Society, Providence, RI, revised edition, 2001.

\bibitem{MillerSturmfels}
E.~Miller and B.~Sturmfels.
\newblock {\em Combinatorial commutative algebra}, volume 227 of {\em Graduate
  Texts in Mathematics}.
\newblock Springer-Verlag, New York, 2005.

\bibitem{Nakamura}
I.~Nakamura.
\newblock Hilbert schemes of abelian group orbits.
\newblock {\em J. Algebraic Geom.}, 10(4):757--779, 2001.

\bibitem{RRZ}
Z.~Reichstein, D.~Rogalski, and J.~J. Zhang.
\newblock Projectively simple rings.
\newblock {\em Adv. Math.}, 203(2):365--407, 2006.

\bibitem{Reid:canonical}
M.~Reid.
\newblock Canonical {$3$}-folds.
\newblock In {\em Journ\'ees de {G}\'eometrie {A}lg\'ebrique d'{A}ngers,
  {J}uillet 1979/{A}lgebraic {G}eometry, {A}ngers, 1979}, pages 273--310.
  Sijthoff \& Noordhoff, Alphen aan den Rijn, 1980.

\bibitem{Reid:young}
M.~Reid.
\newblock Young person's guide to canonical singularities.
\newblock In {\em Algebraic geometry, {B}owdoin, 1985 ({B}runswick, {M}aine,
  1985)}, volume~46 of {\em Proc. Sympos. Pure Math.}, pages 345--414. Amer.
  Math. Soc., Providence, RI, 1987.

\bibitem{SmithStaniszkis}
S.~P. Smith and J.~M. Staniszkis.
\newblock Irreducible representations of the {$4$}-dimensional {S}klyanin
  algebra at points of infinite order.
\newblock {\em J. Algebra}, 160(1):57--86, 1993.

\bibitem{SmithTate}
S.~P. Smith and J.~Tate.
\newblock The center of the {$3$}-dimensional and {$4$}-dimensional {S}klyanin
  algebras.
\newblock {\em K-Theory}, 8(1):19--63, 1994.

\bibitem{StaffordVandenbergh}
J.~T. Stafford and M.~van~den Bergh.
\newblock Noncommutative curves and noncommutative surfaces.
\newblock {\em Bull. Amer. Math. Soc. (N.S.)}, 38(2):171--216 (electronic),
  2001.

\bibitem{Walton:Sdeg}
C.~Walton.
\newblock Degenerate {S}klyanin algebras and generalized twisted homogeneous
  coordinate rings.
\newblock {\em J. Algebra}, 322(7):2508--2527, 2009.

\bibitem{Walton:thesis}
C.~Walton.
\newblock {\em On degenerations and deformations of {S}klyanin algebras}.
\newblock PhD thesis, University of Michigan, deepblue.lib.umich.edu, 2011.

\end{thebibliography}

\end{document}